\documentclass[a4paper, BCOR=0.0mm, DIV=calc]{scrartcl}
\usepackage[T1]{fontenc}
\usepackage[utf8]{inputenc}
\usepackage[ngerman]{babel}
\usepackage[osf,sc]{mathpazo}
\usepackage{textcomp}
\usepackage{microtype}
\usepackage{amsmath, amssymb, amsfonts}
\KOMAoptions{twoside=false, twocolumn=false, headinclude=false, footinclude=false, mpinclude=false, pagesize=auto}
\recalctypearea
\newcommand{\ebinom}[2]{\left(\frac{#1}{#2} \right)}
\begin{document}
\boldmath
\title{Verschiedene Betrachtungen über hypergeometrische Reihen\footnote{
Originaltitel: "`Variae considerationes circa series hypergeometricas"', erstmals publiziert in "`\textit{Nova Acta Academiae Scientarum Imperialis Petropolitinae} 8, 1794, pp. 3-14"', Nachdruck in "`\textit{Opera Omnia}: Series 1, Volume 16, pp. 178 - 192"', Eneström-Nummer E661, übersetzt von: Alexander Aycock, Textsatz: Artur Diener,  im Rahmen des Projektes "`Eulerkreis Mainz"' }}
\unboldmath
\author{Leonhard Euler}
\maketitle
\paragraph{§1}
Nachdem dieses ins Unendliche laufende Produkt
\[
	\frac{P}{Q} = \frac{a(a+2b)}{(a+b)(a+b)} \cdot \frac{(a+2b)(a+4b)}{(a+3b)(a+3b)} \cdot \frac{(a+4b)(a+6b)}{(a+5b)(a+5b)} \cdot \frac{(a+6b)(a+8b)}{(a+7b)(a+7b)}\cdot \mathrm{etc.} 
\]
vorgelegt wurde, ist bekannt, dass
\[
	P = \int\frac{x^{a+b-1}\partial x}{\sqrt{1-x^{2b}}} ~~ \mathrm{und} ~~ Q = \int\frac{x^{a-1}\partial x}{\sqrt{1-x^{2b}}}
\]
ist, wobei diese Integrale von $x=0$ bis $x=1$ erstreckt wurden; dort bemerke man, dass das Glied, welches dem Index $i$ entspricht,
\[
	\frac{(a+(2i-2)b)(a+2ib)}{(a+(2i-1)b)(a+(2i-1)b)}
\]
ist.

\paragraph{§2}
Schon am Beispiel dieses Produktes wollen wir das folgende unbestimmte Produkt betrachten, in welchem die Anzahl der Faktoren gleich $n$ ist, und man setze
\[
	\Delta : n = a(a+2b)(a+4b)(a+6b) \cdot (a+(2n-2)b),
\]
weil ja dieses Produkt wegen der gegebenen Zahlen $a$ und $b$ als eine bestimmte Funktion von $n$ betrachtet werden kann; aus ihrer Gestalt ist daher klar, dass
\[
	\Delta : (n+1) = \Delta : n \cdot (a+2nb)
\]
sein wird und auf ähnliche Weise
\[
	\Delta : (n+2) = \Delta : (n+1) \cdot (a+(2n+2)b)
\]
und so weiter. Daher wird, wenn $i$ eine unendliche große Zahl bezeichnet,
\[
	\Delta : i = a(a+2b)(a+4b)(a+6b)\cdots (a+(2i-2)b)
\]
sein, woher man in gleicher Weise berechnet, dass
\begin{align*}
\Delta : (i+1) &= \Delta : i(a+2ib) \\
\Delta : (i+2) &= \Delta : i(a+2ib)(a+(2i+2)b) \\
\Delta : (i+3) &= \Delta : i(a+2ib)(a+(2i+2)b)(a+(2i+4)b) \\
\mathrm{etc.} &
\end{align*}
sein wird, wo die darüber hinaus hinzukommenden Faktoren als untereinander gleich angesehen werden können; deshalb kann im Allgemeinen
\[
	\Delta : (i+n) = \Delta : i(a+2ib)^n
\]
gesetzt werden, wo, weil $(a+2ib)$ der unmittelbar folgende Faktor ist, mit selbiger Rechtfertigung ein beliebiger der folgenden hätte genommen werden können, woraus wir noch allgemeiner
\[
	\Delta : (i+n) = (\alpha + 2ib)^n \Delta : i
\]
setzen können, während $\alpha$ eine beliebige endliche Zahl bezeichnet, die natürlich in Bezug auf $2ib$ verschwindet.

\paragraph{§3}
Wir wollen nun den Fall des unendlichen Produktes ins Kalkül ziehen, in welchem $n=\frac{1}{2}$ ist, und wir wollen $\Delta : \frac{1}{2} = k$ nennen, welcher Wert sich mit Hilfe der Methode der Interpolation in der Tat immer angeben lässt. Daher wird also durch das Obere
\begin{align*}
\Delta : (1 + \tfrac{1}{2}) &= k(a+b) \\
\Delta : (2 + \tfrac{1}{2}) &= k(a+b)(a+3b) \\
\Delta : (3 + \tfrac{1}{2}) &= k(a+b)(a+3b)(a+5b) \\
\mathrm{etc.}
\end{align*}
sein, woher, indem man ins Unendliche fortschreitet,
\[
	\Delta : (i+\tfrac{1}{2}) = k(a+b)(a+3b)(a+5b) \cdots (a+(2i-1)b) 
\]
sein wird.

\paragraph{§4}
Weil wir also oben schon die Formel für $\Delta : (i+n)$ gegeben haben, werden wir nun für $n=\frac{1}{2}$ gesetzt auch
\[
	\Delta : (i + \tfrac{1}{2}) = \Delta : i \sqrt{\alpha + 2ib}
\]
haben; und so haben wir für dieselbe Formel $\Delta : (i + \frac{1}{2})$ zwei verschiedene Ausdrücke erhalten, und aus diesen folgert man diese Gleichung:
\[
	\Delta : i \sqrt{\alpha + 2ib} = k(a+b)(a+3b)(a+5b) \cdots (a+(2i-1)b)
\]
und daher können wir den Wert des unendlichen Produktes
\[
	(a+b)(a+3b)(a+5b) \cdots (a+(2i-1)b) = \frac{\Delta :i \sqrt{\alpha + 2ib}}{k}
\]
selbst schließen und es wird die Relation zwischen diesem Produkt und dem bekannt, was wir oben duch $\Delta : i$ ausgedrückt haben. Hier ist aber natürlich zu bemerken, dass die Faktoren diese Produktes diese selbst sind, welche den Nenner des eingangs vorgelegten Produktes ergeben; deshalb können wir den Zähler jenes Produktes wie den Nenner durch die gerade gefundenen Werte
\[
	\Delta : i ~~ \mathrm{und} ~~ \frac{\Delta : i\sqrt{\alpha + 2ib}}{k}
\]
ausdrücken.

\paragraph{§5}
Der Zähler des vorgelegten Produktes kann ins Unendliche ausgedehnt so dargestellt werden:
\[
	a(a+2b)^2 (a+4b)^2 \cdots (a+(2i-2)b)^2 (a+2ib)
\]
wo der erste und der letzte Faktor einzeln vorkommen, die übrigen aber alle quadriert. Weil daher
\[
	(\Delta : i)^2 = (a)^2 (a+2b)^2 (a+4b)^2 (a+6b)^2 \cdots (a+(2i-2)b)^2
\]
ist, ist klar, dass jener Zähler $\frac{(\Delta : i)^2}{a}(a+2ib)$ ist. Für den Nenner aber ist per se klar, dass er gleich dem Quadrat des anderen Produktes $(a+b)(a+3b)\mathrm{etc.} $ ist; weil sein Wert mit
\[
	\frac{\Delta : i\sqrt{\alpha + 2ib}}{k}
\]
gefunden wurde, wird der Nenner
\[
	\frac{(\Delta : i)^2(\alpha + 2ib)}{kk}
\]
sein; deshalb haben wir nach Einsetzen dieser Werte für den oben ausgelegten Bruch $\frac{P}{Q}$ diese Gleichung erhalten:
\[
\frac{P}{Q} = \frac{\frac{(\Delta : i)^2 (a + 2ib)}{a}}{\frac{(\Delta : i)^2 (\alpha + 2ib)}{kk}} = \frac{kk(a+2ib)}{a(\alpha + 2ib)} = \frac{kk}{a}
\]
Aus dieser Gleichung wird daher sofort der wahre Wert der interpolierten Formel $k=\Delta : \frac{1}{2}$ bekannt, weil ja
\[
	\Delta : \tfrac{1}{2} = \sqrt{\frac{aP}{Q}}
\]
sein wird und daher weiter die folgenden:
\begin{align*}
	\Delta : (1+\tfrac{1}{2}) &= (a+b)\sqrt{\frac{aP}{Q}} \\
	\Delta : (2+\tfrac{1}{2}) &= (a+b)(a+3b)\sqrt{\frac{aP}{Q}} \\
	\Delta : (3+\tfrac{1}{2}) &= (a+b)(a+3b)(a+5b)\sqrt{\frac{aP}{Q}} \\
	\mathrm{etc.} &
\end{align*}
und diese Interpolation ist umso bemerkenswerter, weil sie ohne Approximation sofort den wahren Wert diese interpolierten Terme liefert.

\paragraph{§6}
Wenn wir darüber hinaus dieses unendliche Produkt, in dem die Faktoren verbunden werden, betrachten und 
\[
	a(a+b)(a+2b)(a+3b) \cdots (a+(i-1)b) = \Gamma : i
\]
setzen, wird
\[
	\Gamma : 2i = a(a+b)(a+2b)(a+3b) \cdots (a+(2i-1)b)
\]
sein, was natürlich das Produkt aus den beiden oberen ist, sodass
\[
	\Gamma : 2i = \frac{(\Delta :i)^n \sqrt{\alpha + 2ib}}{k} 
\]
ist; wenn wir daher die Form $\Gamma : 2i$ benutzen wollten, können wir die Werte der beiden vorhergehenden daraus angeben, weil
\[
	\Delta : i = \frac{k\cdot \Gamma : 2i}{\sqrt{\alpha + 2ib}}
\]
ist, welches der Wert selbst des ersten Produktes
\[
	a(a+2b)(a+4b)(a+6b) ~~ \mathrm{etc.}
\]
ist; der Wert des anderen Produktes
\[
	(a+b)(a+3b)(a+5b) ~~ \mathrm{etc.}
\]
wird
\[
	\frac{\sqrt{\Gamma : 2i \sqrt{\alpha + 2ib}}}{k}
\]
sein.

\paragraph{§7}
Bis hierher haben wir also drei ins Unendliche laufende und miteinander verbundene Produkte betrachtet, die, weil wir sie ja schon gründlicher untersucht haben, hier nochmal vor Augen führen wollen:
\begin{align*}
\mathrm{I.}~~~ & a(a+b)(a+2b)(a+3b) \cdots (a+(i-1)b) = \Gamma : i \\
\mathrm{II.}~~ & a(a+2b)(a+4b)(a+6b) \cdots (a+(2i-2)b) = \Delta : i \\
\mathrm{III.}~ & (a+b)(a+3b)(a+5b) \cdots (a+(2i-1)b) = \Theta : i
\end{align*}
und wird haben gefunden, dass
\[
	\Theta : i = \frac{\Delta : i \sqrt{\alpha + 2ib}}{k}
\]
ist; dann haben wir aber $\Delta : i$ wie $\Theta : i$ auf folgende Weise durch die Funktion $\Gamma : 2i$ ausgedrückt:
\[
	\Delta : i = \sqrt{\frac{k \cdot \Gamma :2i}{\sqrt{\alpha + 2ib}}} ~~ \mathrm{und} ~~ \Theta : i = \sqrt{\frac{\Gamma : 2i\sqrt{\alpha + 2ib}}{k}},
\]
weil ja klar ist, dass
\[
	\Gamma : 2i = \Delta : i \cdot \Theta : i
\]
ist, wo man sich erinnern sollte, dass $k=\Delta : \frac{1}{2}$ ist, was natürlich aus der zweiten Form bestimmt werden muss, indem man die Reihe
\[
	a,~ a(a+2b),~ a(a+2b)(a+4b),~ a(a+2b)(a+4b)(a+6b),~ \mathrm{etc.}
\]
betrachtet, deren Term, der dem Index $\frac{1}{2}$ entspricht, wir mit dem Buchstaben $k$ bezeichnet haben.

\paragraph{§8}
Gleich wollen wir auf diese Formen die allgemeine Methode Progressionen jeder Art durch ihren allgemeinen Term zu summieren anwenden, die sich so verhält, dass nachdem eine beliebige Reihe $A,\, B,\, C,\, D,\, E,\, \mathrm{etc.}$ vorgelegt wurde, deren Term, der dem unbestimmten Index $x$ entspricht, gleich $X$ sei, ihre Summe
\[
	A + B + C + D + \cdots + X,
\]
die wir gleich $S$ nennen wollen,
\[
	S = \int X \partial x + \frac{1}{2}X + \frac{1}{1\cdot 2 \cdot 3}\frac{1}{2}\frac{\partial X}{\partial x} - \frac{1}{1\cdot 2 \cdot 3 \cdot 4 \cdot 5} \frac{1}{6}\frac{\partial^3 X}{\partial x^3} + \frac{1}{1\cdot 2 \cdots 7}\frac{1}{6}\frac{\partial^5 X}{\partial x^5} - \mathrm{etc.}
\]
ist, wo die Brüche $\frac{1}{2},\, \frac{1}{6},\, \frac{1}{6},\, \frac{3}{10},\, \frac{5}{6},\, \mathrm{etc.}$ die Bernoulli-Zahlen sind.

\boldmath
\section*{Entwicklung der ersten Form $(a+b)(a+2b)(a+3b)\cdots (a+(i-1)b) = \Gamma : i$}
\unboldmath
\paragraph{§9}
Weil hier die Anzahl der Faktoren als unendlich betrachtet wird, wollen wir, damit wir die Summierungsmethode auf sie anwenden können, dieselbe Form, die aus einer endlichen Anzahl an Termen, z.\,B. gleich x, betrachten und wollen auf ähnliche Weise
\[
	a(a+b)(a+2b)(a+3b) \cdots (a+(x-1)b) = \Gamma : x
\]
setzen. Nun aber wollen wir, damit wir anstelle dieses Produktes eine summierende Reihe erhalten, Logarithmen nehmen und es wird
\[
	\log{\Gamma : x} = \log{a} + \log{(a+b)} + \log{(a+2b)} + \log{(a+3b)} + \cdots + \log{(a+(x-1)b)}
\]
sein; nachdem also die Summe dieser (Form) gefunden wurde, wird sie den Logarithmus der Formel $\Gamma : x$ geben und daher die Formel $\Gamma : x$ selbst; wenn in ihr $x=i$ gesetzt wird, wird man die Formel $\Gamma : i$ erhalten, welchen Wert wir hauptsächlich in den oberen Paragraphen betrachtet haben. Daher wird also, nachdem ein Vergleich mit der allgemeinen Reihe aufgestellt wurde, 
\[
	X = \log{(a+(x-1)b)}
\]
sein und die Summe selbst
\[
	S = \log{\Gamma : x},
\]
oder es wird
\[
	X = \log{(a-b+bx)}
\]
sein, woher man
\[
	\int X \partial x = \int \partial x \log{(a-b+bx)}
\]
berechnet.

\paragraph{§10}
Weil also
\[
	\int \partial z \log{z} = z\log{z} - z
\]
und
\[
	\int \partial y \log{(a+y)} = (a+y)\log{(a+y)} - (a+y)
\]
ist, wird nun, indem man anstelle von $y$ gleich $bx$ schreibt,
\[
	\int b \partial x \log{(a+bx)} = (a+bx)\log{(a+bx)} - a - bx
\]
sein und daher
\[
	\int \partial x \log{(a+bx)} = \frac{a+bx}{b}\log{(a+bx)} - \frac{a}{b} - x,
\]
woher man für unseren Fall berechnet, dass
\[
	\int X \partial x = \frac{(a-b+bx)}{b}\log{(a-b+bx)} - \frac{a}{b} + 1 - x
\]
sein wird, wo sich im letzten Teil das konstante Glied $\frac{a}{b} - 1$ weglassen lässt, weil dieser Ausdruck per se eine Konstante unbestimmter Größe erfordert, die man schließlich aus der Gestalt der Reihe bestimmen muss. Darauf wird aber
\[
	\frac{\partial X}{\partial x} = \frac{b}{a-b+bx}
\]
sein, dann aber weiter
\[
	\frac{\partial^3 X}{\partial x^3} = \frac{2b^3}{(a-b+bx)^3}, ~~ \frac{\partial^5 X}{\partial x^5} = \frac{2\cdot 3 \cdot 4b^5}{(a-b+bx)^5}, ~~ \mathrm{etc,}
\]
nach Benutzung welcher Werte
\begin{align*}
	2\Gamma : x =& A + \left( \frac{a}{b} - \frac{1}{2} + x\right) \log{(a-b+bx)} - x + \frac{1}{1\cdot 2 \cdot 3}\cdot \frac{1}{2} \cdot \frac{b}{a-b+bx} \\
	&-\frac{1}{3\cdot 4\cdot 5}\cdot \frac{1}{6} \cdot \frac{b^3}{(a-b+bx)^3} + \frac{1}{5\cdot 6 \cdot 7}\cdot \frac{1}{6}\cdot \frac{b^5}{(a-b+bx)^5} \\
	&-\frac{1}{7\cdot 8\cdot 9} \cdot \frac{3}{10} \cdot \frac{b^7}{(a-b+bx)^7} + \frac{1}{9\cdot 10 \cdot 11} \cdot \frac{5}{6} \cdot \frac{b^9}{(a-b+bx)^9} - \mathrm{etc,}
\end{align*}
wo der Buchstabe $A$ eine aus der Gestalt der Reihe zu bestimmende Konstante bezeichnet.

\paragraph{§11}
Diese Konstante $A$ muss aber aus einem Fall, in welchem die Summe bekannt ist, bestimmt werden, was also aus dem Fall $x=0$ geschehen könnte, in dem ja die Summe gleich $0$ sein muss; aber wäre also daher
\begin{align*}
A =& \left( \frac{a}{b} - \frac{1}{2} \right)\log{(a-b)} + \frac{1}{1\cdot 2 \cdot 3} \cdot \frac{1}{2}\cdot \frac{b}{a-b} - \frac{1}{3\cdot 4 \cdot 5}\cdot \frac{1}{6} \cdot \frac{b^3}{(a-b)^3} \\
&+ \frac{1}{5\cdot 6 \cdot 7}\cdot \frac{1}{6}\cdot \frac{b^5}{(a-b)^5} - \mathrm{etc.}
\end{align*}
Weil aber diese Reihe kaum kovergiert und sogar im Fall $b=a$ alle Terme unendlich werden würden, lässt sich daher natürlich nichts gewinnen. Wenn wir aber $x=1$ nehmen wollten, müsste die Summe $\log{a}$ hervorgehen, woher sich in gleicher Weise kaum irgendetwas für unsere Aufgabe schließen ließe, weil man immer zu einer unendlichen Reihe gelangen würde, deren Summe man erst erforschen müsste, bei welcher Aufgabe zwar vielleicht das, was ich einst pber Reihen, die Bernoulli-Zahlen involvieren, bemerkt habe, einen gewissen Nutzen leisten könnte, was aber wiederum mit nicht gerade wenig Arbeit verbunden wäre.

\paragraph{§12}
Weil wir nämlich bei der gegenwärtigen Aufgabe hauptsächlich auf den Wert $\Gamma : i$ schauen, wird es genügen, dass sofort anstelle von $x$ eine unendliche Zahl gesetzt wird. Es sei daher also $x=i$, während $i$ eine unendlich große Zahl bezeichnet, und unsere Gleichung nimmt diese Form an:
\[
	\log{\Gamma : i} = A + \left( \frac{a}{b} - \frac{1}{2} + i\right) \log{(a-b+bi)} - i,
\]
woher die Konstante $A$ von selbst bestimmt wird, welche wir deshalb als schon bekannt ansehen werden. Daher werden wir, indem wir zu Zahlen zurückgehen, wo wir anstelle von $A$ aber $\log{A}$ geschrieben verstehen wollen, zu diesem Ausdruck gelangen
\[
	\Gamma : i = A(a-b+bi)^{\frac{a}{b} - \frac{1}{2} + i}e^{-i}.
\]
Hier sollte freilich die Potenz des Exponenten $i$ getrennt dargestellt werden, und zwar auf diese Weise:
\[
	\Gamma : i = A(a-b+bi)^{\frac{a}{b} - \frac{1}{2}}(a-b+bi)^{i}e^{-i}.
\]

\section*{Entwicklung der beiden übrigen Formeln}
\paragraph{§13}
Die zweite Form unterscheidet sich von der ersten nur darin, dass anstelle von $b$ hier $2b$ geschrieben werden muss, woher wir auf eine neue Entwicklung verzichten können; aber wir wollen  anstelle der Konstante $A$ hier $B$ schreiben, weil ja noch nicht bekannt ist, wie der Buchstabe $b$ in die Konstante $A$ eingeht. Auf diese Weise werden wir also sofort
\[
	\Delta : i = B(a-2b+2bi)^{\frac{a}{2b} - \frac{1}{2}}(a-2b+2bi)^i e^{-1}
\]
haben. Auf ähnliche Weise ist klar, dass aus dieser zweiten Form die dritte entsteht, wenn nur anstelle von $a$ gleich $a+b$ geschrieben wird, woher wir, indem wir anstelle von $B$ die Konstante $C$ einführen,
\[
	\Theta : i =C(a-b+2bi)^{\frac{a}{2b}}(a-b+2bi)^i e^{-i}
\]
haben werden. Dort bemerke man, dass der Buchstabe $e$ hier für die Zahl gesetzt wurde, deren hyperbolischer Logarithmus gleich $1$ ist.

\section*{Daher entstehende Schlussfolgerungen}
\paragraph{§14}
Wir werden nun sehen, wie diese neuen Bestimmungen sich in Anbetracht der oben gefundenen Relationen verhalten werden; weil daher aus diesen neuen Werten
\[
	\Gamma : 2i = A(a-b+2bi)^{\frac{a}{b} - \frac{1}{2}}(a-b+2bi)^{2i}e^{-2i}
\]
ist, werden wir, weil wir
\[
	\Gamma : 2i = \Delta : i \cdot \Theta : i
\]
gefunden haben, wenn wir hier überall die gerade gefundenen Werte einsetzen, für diese Gleichung zuerst das Produkt
\[
	\Delta : i \cdot \Theta : i = B C (a-2b+2bi)^{\frac{a}{2b}-\frac{1}{2}}(a-b+2bi)^{\frac{a}{2b}}(a-b+2bi)^i(a-b+2bi)^i e^{-2i}
\]
haben; weil dieses (Produkt) jenem Wert $\Gamma : 2i$ gleich sein muss, wird, wenn durch die Faktoren, welche es gemeinsam hat, auf beiden Seiten geteilt wird, diese Gleichung
\[
	A(a-b+2bi)^{\frac{a}{2b}-\frac{1}{2}}(a-b+2bi)^i = B C (a-2b+2bi)^{\frac{a}{2b} - \frac{1}{2}}(a-2b+2bi)^i
\]
hervorgehen.

\paragraph{§15}
Wir wollen nun diese Gleichung auf beiden Seiten durch $(a-2b+2bi)^i$ teilen, und weil
\[
	\frac{a-b+2bi}{a-2b+2bi} = 1 + \frac{b}{a-2b+2bi} = 1 + \frac{1}{2i}
\]
gleich $i$ ist, wird durch eine gewöhnliche Auflösung
\[
	\left( 1 + \frac{1}{2i}\right)^i = e^{\frac{1}{2}}
\]
sein, woher unsere Gleichung auf diese Form
\[
	A(a-b+2bi)^{\frac{a}{2b} - \frac{1}{2}}e^{\frac{1}{2}} = BC(a-2b+2bi)^{\frac{a}{2b}-\frac{1}{2}}
\]
zurückgeführt wird, wo sich auch die letzten Faktoren aufheben, weil
\[
	\ebinom{a-b+2bi}{a-2b+2bi}^{\frac{a}{2b} - \frac{1}{2}} = \left( 1 + \frac{1}{2i}\right)^{\frac{a}{2b}-\frac{1}{2}} = 1
\]
ist, so dass man zu dieser einfachen Gleichheit gelangt:
\[
	Ae^{\frac{1}{2}} = B C.
\]

\paragraph{§16}
Weil wir darauf oben gefunden haben, dass
\[
	\Theta : i = \frac{\Delta : i \sqrt{\alpha + 2ib}}{k}
\]
ist oder
\[
	\frac{\Theta : i}{\Delta : i} = \frac{\sqrt{\alpha + 2ib}}{k},
\]
wollen wir den für $\Theta : i$ gefundenen Wert durch $\Delta : i$ teilen und wir werden
\[
	\frac{\Theta : i}{\Delta : i} = \frac{C}{B}\sqrt{a-2b+2bi}\ebinom{a-b+2bi}{a-2b+2bi}^i = \frac{C}{B}\sqrt{e(a-2b+2bi)}
\]
finden. Es wird also
\[
	\frac{\sqrt{\alpha + 2ib}}{k} = \frac{C}{B}\sqrt{e(a-2b+2bi)}
\]
sein, oder
\[
	\frac{1}{k} = \frac{C}{B}\sqrt{\frac{e(a-2b+2bi)}{\alpha + 2ib}} = \frac{C}{B}\sqrt{e}
\]
oder es wird
\[
	B = Ck\sqrt{e}
\]
sein.

\paragraph{§17}
Wir haben also zwei Relationen solcher Art zwischen jenen drei Konstanten $A, B, C$ erhalten, sodass, wenn eine einzige bekannt wäre, aus ihr die beiden übrigen bestimmt werden könnten. Weil nämlich
\[
	A = \frac{BC}{\sqrt{e}} ~~ \mathrm{und} ~~ B = Ck\sqrt{e}
\]
ist, werden, wenn wir die Konstante $A$ als schon bekannt betrachten, die beiden übrigen auf folgende Weise bestimmt werden. Weil $B = CK\sqrt{e}$ ist, gibt dieser Wert in die erste Gleichung eingesetzt $CCk=A$, woher man $C = \sqrt{\frac{A}{k}}$ findet und daher weiter $B = \sqrt{k}Ae$. Trotzdem ist daher nicht klar, auf welche Weise die Konstante $A$ uneingeschränkt bestimmt werden kann, und daher wird zu jener Summation der logarithmischen Reihe selbst zurückzukehren sein, die wir oben mit dem Buchstaben $A$ bezeichnet haben, wo aber anstelle von $A$ gleich $\log{A}$ zu schreiben sein wird. Und daher haben wir nur gewonnen, dass, wenn die beiden übrigen Formen auf ähnliche Weise durch logarithmische Reihen entwickelt werden, die dort hinzuzufügenden Konstanten, natürlich $\log{B}$ und $\log{C}$, zugleich bekannt werden.

\paragraph{§18}
Es ist übrig, dass wir noch ein paar Dinge über die Werte des Buchstaben $k$ hinzufügen, an den wir schon oben erinnert haben, dass er durch Interpolation gefunden werden muss. Dennoch kann dieser Buchstabe auch aus dem Vergleich selbst der Formeln $\Delta : i$ und $\Theta : i$ uneingeschränkt durch bestimmt Quadraturen bestimmt werden. Weil nämlich
\[
	k = \frac{\Delta : i}{\Theta : i}\sqrt{\alpha + 2ib}
\]
ist und daher
\[
	kk = \frac{(\Delta : i)^2 (\alpha + 2ib)}{(\Theta : i)^2},
\]
geht hier, wenn wir anstelle von $\Delta : i$ und $\Theta : i$ die unendlichen Produkte selbst einsetzen, weil ja beide aus $i$ Faktoren bestehen, und wenn wir den Faktor des Zählers getrennt ausdrücken, aber im Zähler ein einziger Faktor $\alpha + 2ib$ darüber hinaus ein; nachdem das festgesetzt wurde, werden wir zu dem folgenden bestimmten Produkt gelangen:
\[
	kk = a \cdot \frac{a(a+2b)(a+2b)(a+4b)(a+4b)(a+6b)}{(a+b)(a+b)(a+3b)(a+3b)(a+5b)(a+5b)} ~~\mathrm{etc.}
\]

\paragraph{§19}
Damit wir aber den wahren Wert diese unendlichen Produktes finden, ist zu bemerken, wenn die Buchstaben $P$ und $Q$ die folgenden Integralformeln bezeichnen
\[
	P = \int\frac{x^{p-1}\partial x}{(1-x^n)^{1-\frac{m}{n}}} ~~ \mathrm{und} ~~ Q = \int\frac{x^{q-1}\partial x}{(1-x^n)^{1-\frac{m}{n}}},
\]
welche Integrale zu verstehen sind, dass sie von $x=0$ bis $x=1$ erstreckt werden, dass dann durch ein unendliches Produkt 
\[
	\frac{P}{Q} = \frac{q(m+p)}{p(m+q)} \cdot \frac{(q+n)(m+p+n)}{(p+n)(m+q+n)} \cdot \frac{(q+2n)(m+p+2n)}{(p+2n)(m+q+2n)} \cdot \mathrm{etc.}
\]
sein wird, welches Produkt leicht auf unsere Form zurückgeführt wird, indem man
\[
	q=a,~ p=a+b,~ m=b,~ n=2b
\]
nimmt, sodass für unseren Fall
\[
	P = \int\frac{x^{a+b-1}\partial x}{\sqrt{1-x^{2b}}} ~~ \mathrm{und} ~~ Q = \int\frac{x^{a-1}\partial x}{\sqrt{1-x^{2b}}}
\]
wird; dann wird aber
\[
	kk = \frac{aP}{Q}
\]
sein und daher
\[
	k = \sqrt{\frac{aP}{Q}}
\]
und so haben wir denselben Wert $k$ auf andere Art gefunden, den wir oben schon angemerkt haben.

\paragraph{§20}
Sowie aber $k = \Delta : \frac{1}{2}$ ist, können wir auf ähnliche Weise für die beiden übrigen Formen die Werte $\Gamma : \frac{1}{2}$ und $\Theta : \frac{1}{2}$ angeben. Weil nämlich die Form $\Gamma$ aus der Form $\Delta$ ensteht, wenn in dieser anstelle von $b$ gleich $\frac{1}{2}b$ geschrieben wird, die Form $\Theta$ aber aus $\Delta$ entsteht, wenn anstelle von $a$ gleich $a+b$ geschrieben wird, wird, nachdem das bemerkt wurde,
\[
	\Gamma : \frac{1}{2} = \sqrt{a\frac{\int \frac{x^{a+\frac{1}{2}b-1}\partial x}{\sqrt{1-x^b}}}{\int\frac{x^{a-1}\partial x}{\sqrt{1-x^{b}}}}}
\]
und
\[
	\Theta : \frac{1}{2} = \sqrt{(a+b)\frac{\int\frac{x^{a+2b-1}\partial x}{\sqrt{1-x^{2b}}}}{\int\frac{x^{a+b-1}\partial x}{\sqrt{1-x^{2b}}}}}
\]
sein. Man sieht aber leicht an, dass der Wert $\Theta : \frac{1}{2}$ genauso in unsere Rechnung eingeführt werden könnte wie $\Delta : \frac{1}{2} = k$, weil
\[
	\Delta : \frac{1}{2} \cdot \Theta : \frac{1}{2} = a
\]
ist. Nachdem nämlich jene Integralwerte miteinander mutlipliziert wurden, geht
\[
	\Delta : \frac{1}{2} \cdot \Theta : \frac{1}{2} = \sqrt{\frac{a(a+b)\int\frac{x^{a+2b-1}\partial x}{\sqrt{1-x^{2b}}}}{\int\frac{x^{a-1}\partial x}{(1-x^{2b})}}}
\]
hervor; aus einer allbekannten Reduktion solcher Integrale ist aber bekannt, dass
\[
	\int \frac{x^{a+2b-1}\partial x}{\sqrt{1-x^{2b}}} = \frac{a}{a+b}\int\frac{x^{a-1}\partial x}{\sqrt{1-x^{2b}}}
\]
ist, natürlich für die Integrationsgrenzen $x=0$ und $x=1$, und so ist klar, dass
\[
\Delta : \frac{1}{2} \cdot \Theta : \frac{1}{2} = a
\]
sein wird. Wie sich aber der Wert $\Gamma : \frac{1}{2}$ auf beide übrigen bezieht, kann auf keine Weise bestimmt werden.
\end{document}